\input amstex
\documentstyle{amsppt}
\magnification=\magstep1

\input amssym.tex
\input amssym.def
\input epsf.sty

\hsize=6.5truein
\vsize=8.9truein

\NoBlackBoxes


\def\bs{{\backslash}}
\def\ds{{\displaystyle}}
\def\dscdot{{\ds\cdot}}

\def\det{{\text {\rm det}}}

\def\a {{\alpha}}
\def\b {{\beta}}
\def\d {{\delta}}
\def\D {{\Delta}}
\def\g {{\gamma}}
\def\la {{\lambda}}

\def\p {{\partial}}

\def\Si {{\Sigma}}
\def\vp{{\varphi}}

\def\BC{{\Bbb C}}
\def\BK{{\Bbb K}}
\def\BN{{\Bbb N}}
\def\BP{{\Bbb P}}

\def\BR{{\Bbb R}}

\def\cB{{\Cal B}}

\def\cO{{\Cal O}}
\def\cP{{\Cal P}}
\def\cT{{\Cal T}}

\def\hDelta{{\widehat\Delta}}
\def\hf {{\widehat f}}
\def\hcO{{\widehat\cO}}
\def\hvp{{\widehat\vp}}

\def\oU{{\overline U}}

\def\ua{{\underline a}}
\def\ux{{\underline x}}
\def\uvp{{\underline\vp}}
\def\uE{{\underline E}}

\document
\topmatter

\title Control of radii of convergence and extension
of subanalytic functions \endtitle
\leftheadtext{\sevenrm Edward Bierstone}
\rightheadtext{\sevenrm Control of radii of convergence and
extension of subanalytic functions}

\author
Edward Bierstone
\endauthor

\address
{Department of Mathematics, University of Toronto,
Toronto, Ontario, Canada M5S 3G3}
\endaddress

\thanks
The author's research was partially supported by NSERC
grant 0GP0009070.
\endthanks


\abstract
Let $g$: $U\to \BR$ denote a real analytic function on an open subset
$U$ of $\BR^n$, and let 
$\Si \subset \p U$ denote the points where $g$ does not admit a 
local analytic extension. We show that if $g$ is semialgebraic 
(respectively, globally subanalytic), then $\Si$ is semialgebraic
(respectively, subanalytic) and $g$ extends to a neighbourhood of
$\oU\bs \Si$ as an analytic function that is semialgebraic 
(respectively, globally subanalytic). (In the general subanalytic case,
$\Si$ is not necessarily subanalytic.) Our proof depends on controlling
the radii of convergence of power series $G$ centred at points $b$ in
the image of an analytic mapping $\vp$, in terms of the radii of convergence
of $G\circ\hvp_a$ at points $a\in\vp^{-1}(b)$, where $\hvp_a$ denotes the
Taylor expansion of $\vp$ at $a$.
\endabstract

\endtopmatter

\baselineskip=18pt

\head 1. Introduction \endhead
Let $g$: $U\to \BR$ denote a real-analytic function defined on
an open subset $U$ of $\BR^n$, such that $g$ is {\it semialgebraic}
(respectively, {\it subanalytic}); i.e., the graph of $g$ is
semialgebraic (respectively, subanalytic) as a subset of $\BR^n\times\BR$.
Chris Miller asked me the following questions about extension of
the domain of $g$:

(1.1) Let $\Si$ denote the subset of the boundary $\p U$ of $U$ consisting
of points where $g$ does not admit an analytic extension (to some
neighbourhood).
Is $\Si$ a closed semialgebraic (respectively, subanalytic) subset
of $\p U$?

(1.2) If so, can $g$ be extended to a neighbourhood of $\oU\bs \Si$
as an analytic function that is semialgebraic (respectively,
subanalytic)?

Note that if $g$: $U\to \BR$ is semialgebraic (respectively, 
subanalytic and bounded),
then $U$ is semialgebraic (respectively, subanalytic).
If $g$ is semialgebraic, then a local analytic extension at a point
of $\p U$ is semialgebraic (on a suitable neighbourhood).
An analytic function that is semialgebraic is called
{\it algebraic} or {\it Nash}.

We will show that the answer to (1.1) is ``yes'' if $g$ is a semialgebraic
function or a bounded (or global) subanalytic function (Theorem 2.3 below) ---
the proof is a direct application of the ``graphic point'' argument
used in \cite{BM, \S7} to show that the set of smooth points of a subanalytic
set is subanalytic.
The answer to (1.1) is ``no'' in the general subanalytic case (Example 2.1).

The answer to (1.2) is ``yes'' under the 
hypotheses above (Theorem 2.4).
Our proof depends on controlling the radii of convergence
of $g$ at points of $\p U\bs \Si$. This essentially means
controlling the radii of
convergence of power series $G$ centred at points $b$ in the image
of an analytic mapping $\vp$, in terms of the radii of convergence
of $G\circ\hvp_a$ at points $a\in\vp^{-1}(b)$.
($G\circ \hvp_a$ denotes the formal composition of $G$ with the
Taylor expansion $\hvp_a$ of $\vp$ at $a$.)
Let $\BK=\BR$ or $\BC$.
Let $V$ be an open neighbourhood of $0$ in $\BK^n$, and let
$\vp$: $V\to\BK^n$ denote an analytic mapping whose Jacobian
determinant
$$
\D := \det {\p (\vp_1,\ldots,\vp_n)\over \p(x_1,\ldots,x_n)}
$$
does not vanish identically in a neighbourhood of $0$.
We write $\cO_b$ to denote the ring of germs of analytic
functions at a point $b$ of $\BK^n$, and $\hcO_b$ to denote
the formal completion of $\cO_b$; i.e., $\hcO_b$ is the ring of
formal power series centred at $b$.
Tougeron \cite{T, 5.10} and Chaumat and Chollet \cite{CC, \S17} have shown
that, for a given point $a\in V$, there exist constants $\la\in\BN$
and $c>0$ such that, if $G\in\hcO_{\vp(a)}$ and $F=G\circ\hvp_a$
converges in a ball $|x-a|<r$, where $r\le 1$, then $G$ has
radius of convergence $r_G\ge cr^\la$.
The power $\la$ can be chosen to be independent of $a$,
but $c$ cannot be chosen uniformly:

\example{Example 1.3}
Let $\vp$: $\BC\to \BC$ denote the mapping $\vp(x)=x^2$.
If $G\in\hcO_0$ and $F=G\circ \hvp_0$ converges in
a disk $|x|<r$, then $G$ converges in $|y|<r^2$.
On the other hand, if $f(x)=x$ and $a\ne 0$, then we can solve
$f=g\circ \vp$ for $g\in \cO_{\vp(a)}$, but $g$ converges only
in $|y-\vp(a)|<|a|^2$.
\endexample

We will show that, in general, the lack of uniformity in
the constant $c$ is of the same nature as in the preceding
elementary example:

\proclaim {Theorem 1.4}
Let $\vp$: $V\to \BK^n$ denote an algebraic (respectively,
analytic) mapping from an open neighbourhood $V$ of $0$ in $\BK^n$.
Assume that the Jacobian determinant $\D$ of $\vp$ does not vanish
identically in a neighbourhood of $0$.
Then we can find a neighbourhood $W$ of $\, 0$ in $V$, constants
$\d>0$ and $\la\in\BN$, and a finite filtration of $W$ by closed
algebraic (respectively, analytic) subsets,
$$
W = X_0 \supset X_1\supset \cdots \supset X_{k+1} = \emptyset ,
$$
such that, for all $j=0, \ldots, k$, there exists $\a=\a(j)\in\BN^n$
with the properties that, for every $a\in X_j\bs X_{j+1}$:

{\rm (1)} $D^\a \D(a)\ne 0$, where $D^\a$ denotes the
partial derivative $\p^{|\a|}/\p x_1^{\a_1}\cdots \p x_n^{\a_n}$, with
$\a = (\a_1,\ldots, \a_n)$, $|\a|=\a_1+\cdots+\a_n$.

{\rm (2)} If $G\in\hcO_{\vp(a)}$ and $F=G\circ \hvp_a$ converges
in a ball $|x-a|<r$, where $r\le 1$, then $G$ has radius
of convergence
$$
r_G \ge \d |D^\a \D(a)|^2 r^\la .
$$
\endproclaim

\proclaim {Corollary 1.5}
Let $U$ be an open subset of $\BR^n$ and let $g$: $U\to\BR$
denote an analytic function that is semialgebraic (respectively,
globally subanalytic).
Let $\Si$ be a closed semialgebraic (respectively, subanalytic)
subset of $\p U$, and assume that $g$ admits a local analytic
extension at each point of $\p U \bs \Si$.
Then there is a partition $\cP=\{S_i\}$ of $\p U\bs\Si$ such that
$\cP$ is finite (respectively, locally finite in $\p U$), 
each $S_i$ is a semialgebraic (respectively, bounded
subanalytic) subset of $\BR^n$, and, for each $i$,
there is a continuous semialgebraic (respectively, subanalytic)
function $\rho_i$: $S_i\to\BR$ such that, for all $b\in S_i$,
$g$ extends to an analytic function that is semialgebraic
(respectively, subanalytic) on the ball of radius $\rho_i(b)$ centred
at $b$.
\endproclaim

We sketch a proof of Theorem 1.4 in Section 3 below;
the theorem follows immediately from estimates used by Augustin
Mouze in \cite{M, \S I.2} to prove a more precise version of
the theorem of Chaumat and Chollet \cite{CC}.
In Section 2, we deduce Corollary 1.5 from Theorem 1.4 and
we answer Miller's questions.
I am happy to acknowledge my discussions with Anne-Marie Chollet,
Augustin Mouze and Vincent Thilliez concerning \cite{CC} and \cite{M},
and to thank Artur Piekosz for pointing out some inaccuracies
in an earlier manuscript.

\head 2. Extension of semialgebraic and subanalytic functions\endhead

The following example shows that the answer to (1.1) is ``no''
in general for a subanalytic function on a subanalytic domain. 

\example{Example 2.1}
Let $U=\{ (x,y)\in\BR^2:\, x>0,\, y>0\}$ and let
$$
g(x,y) = {1\over x^2\left( \sin^2 \displaystyle{{1\over x}} + y^2\right)} .
$$
The denominator of $g(x,y)$ has zero-set
$$
Z=\{x=0\} \cup \left\{ y=0,\ \sin {1\over x} = 0\right\} ,
$$
and is analytic except at points where $x=0$.
Therefore, $g(x,y)$ is subanalytic.
(In fact, it is semianalytic.)
But $\Si= Z\cap \{ x\ge 0,\, y\ge 0\}$ is not subanalytic.
\endexample

Let $g$: $U\to\BR$ denote an analytic function on an open
subset $U$ of $\BR^n$ such that $g$ is semialgebraic
(respectively, subanalytic).
We compactify $\BR^n$ to real projective space $\BP^n$, and
we compactify $\BR$ to $\BP^1= S^1$; write $\infty$ for the point
at infinity of the latter.
The following conditions are obviously equivalent:

(2.2)(i) $g$ is semialgebraic.

(ii) The graph of $g$ is semialgebraic as a subset of $\BR^n\times S^1$.

(iii) The graph of $g$ is semialgebraic as a subset of $\BP^n\times S^1$.

\smallskip

For the subanalytic analogues of these conditions, we have
(iii) $\Rightarrow$ (ii) $\Rightarrow$ (i).
We say that $g$ is a {\it global subanalytic} function if it
satisfies the analogue of (ii); i.e., the graph of $g$ is
subanalytic as a subset of $\BR^n\times S^1$.
Clearly, a bounded subanalytic function is globally subanalytic.

\proclaim{Theorem 2.3}
Let $U$ be an open subset of $\BR^n$ and let $g$: $U\to\BR$
denote an analytic function.
Let $\Si\subset \p U$ denote the set of points at which $g$
does not admit an analytic extension.
If $g$ is semialgebraic (respectively, globally subanalytic),
then $\Si$ is a closed semialgebraic (respectively, subanalytic)
set.
\endproclaim

\demo{Proof}
Let $X$ denote the closure of the graph of $g$ in $\BP^n\times S^1$.
Since the question is local, if $g$ is globally subanalytic,
then we can reduce to the case that the graph of $g$ is subanalytic in
$\BP^n\times S^1$.
So we assume that $X$ is semialgebraic (respectively, subanalytic).

By the uniformization theorem for semialgebraic or subanalytic
sets \cite{BM, Theorem 0.1}, there is a compact real algebraic (respectively,
analytic) manifold $M$ of dimension $n$, and an algebraic (respectively,
analytic) mapping $\Phi=(\vp,f)$: $M\to\BP^n\times S^1$ such that
$\Phi(M)=X$.
We can assume that $\Phi$ (and therefore $\vp$) has maximum
rank $n$ on each component of $M$.
There is a bound $s$ on the number of connected components
of the fibres of $\vp$ \cite{BM, Theorem 3.14}.
Let $M_\vp^s$ denote the {\it $s$-fold fibre-product} of $M$ with
respect to $\vp$; i.e.,
$$
M_\vp^s := \big\{ \ux = (x^1,\ldots,x^s)\in M^s:\ 
\vp(x^1) = \cdots = \vp (x^s)\big\}
$$
and let $\uvp$: $M_\vp^s \to\BP^n$ denote the mapping induced by $\vp$.

We say that $\ua\in M_\vp^s$ is an $s$-{\it fold graphic point}
of $\Phi$ if there exists $g^\ua \in \cO_{\uvp(\ua)}$ such that
$f_{a^i} = g^\ua \circ \vp_{a^i}$, $i=1,\ldots,s$, where $\ua=
(a^1,\ldots, a^s)$ (and $f_{a^i}$, $\vp_{a^i}$ denote the germs
of $f$, $\vp$ at $a^i$).
Let $\uE \subset M_\vp^s$ denote the subset of points
that are {\it not} $s$-fold graphic points of $\Phi$.
Then $\uE$ is a closed algebraic (respectively, analytic)
subset of $M_\vp^s$ \cite{BM, Corollary 7.13}.
But $g$ extends to a neighbourhood of a point $b\in \partial U$
precisely when $b\not\in \uvp(\uE)$ and $(b,\infty)\not\in X$.
$\square$
\enddemo

\demo{Proof of Corollary 1.5}
Let $g$: $U\to\BR$ be 
an analytic function that is semialgebraic (respectively,
globally subanalytic).
Let $Y\subset\BP^n\times S^1$ (respectively,
$Y\subset\BR^n\times S^1$) denote the closure of the graph
of $g$.
By the uniformization theorem \cite{BM, Theorem 0.1}, there is an algebraic
(respectively, analytic) manifold $M$ and a proper algebraic
(respectively, analytic) mapping $\Phi=(\vp,f)$: $M\to \BP^n\times S^1$
(respectively,  $\Phi=(\vp,f)$: $M\to \BR^n\times S^1$)
such that $\Phi(M)=Y$.
Again, we can assume that $\Phi$ (and hence $\vp$) has maximum rank
$n$ on each component of $M$.
Suppose $a\in M$ and $b=\vp(a) \in \p U$.
If $g$ extends to an analytic function in a neighbourhood
of $b$, then $f_a=g_b\circ \vp_a$, where $g_b$ denotes the
germ of an extension at $b$.
Corollary 1.5 is now a simple consequence of Theorem 1.4.
$\square$
\enddemo

\proclaim {Theorem 2.4}
Let $U$ be an open subset of $\BR^n$ and let $g$:
$U\to \BR$ be an analytic function that is semialgebraic
(respectively, globally subanalytic).
Let $\Si$ be a closed semialgebraic (respectively, subanalytic)
subset of $\partial U$.
Assume that $g$ admits a local analytic extension at each point
of $\p U\bs \Si$.
Then $g$ extends to a neighbourhood of $\oU\bs\Si$ as an analytic
function that is semialgebraic (respectively, globally subanalytic).
\endproclaim

\demo{Proof}
Consider a partition $\cP=\{ S_i\}$ of $\p U\bs\Si$
and continuous semialgebraic (respectively, subanalytic)
functions $\rho_i$: $S_i\to \BR$ satisfying the
conclusion of Corollary 1.5.
For each $i$, set
$$
U_i = \bigcup_{b\in S_i} \big\{ y:\ |y-b| < \rho_i (b)\big\} ;
$$
then $U_i$ is an open semialgebraic (respectively, subanalytic)
neighbourhood of $S_i$.
Let $V'=\cup_i U_i$.
Then $V'$ is an open semialgebraic (respectively, subanalytic)
neighbourhood of $\p U\bs \Si$.
By Corollary 1.5, $g$ extends to $U\cup V'$ but, in general,
only as a multivalued function.

Let $\cT$ denote a finite semialgebraic triangulation of $\BP^n$
(respectively, a locally
finite subanalytic triangularization of $\BR^n$) that is compatible
with $U$, $\Si$ and each $U_i$ \cite{H} and let
$\cB(\cT)$ denote the barycentric subdivision of $\cT$.
Let $V$ denote the union of all open simplices of $\cB(\cT)$
that are adherent to $\p U\bs \Si$.
(An {\it open simplex} means a simplex minus its boundary.)
Then $V$ is an open neighbourhood of $\p U\bs \Si$, and $g$
extends to $U\cup V$ as an analytic function that is
semialgebraic (respectively, subanalytic).
$\square$
\enddemo

\remark{Remarks 2.5}
(1) The assertions of Theorems 2.3 and 2.4 under the weaker
hypotheses that $g$ is semialgebraic (respectively, 
globally subanalytic)
but not necessarily analytic throughout $U$ follow immediately
from the theorems as stated.

(2) There are analogues of Theorems 2.3 and 2.4 for a semianalytic
function $g$: $U\to\BR$ that can be proved using a real normalization
of an analytic hypersurface containing the graph of $g$ in a
neighbourhood of a point of its closure. \endremark

\head 3. Control of radii of convergence\endhead
In this section, we show that Theorem 1.4 follows from
the estimates used by Mouze \cite{M, \S I.2} to prove his
generalization of the theorem of Tougeron and
Chaumat-Chollet.
Let $\vp$: $V\to\BK^n$ denote an algebraic
(respectively, analytic) mapping, where $V$ is a
neighbourhood of $0$ in $\BK^n$.
($\BK=\BR$ or $\BC$.)
Let
$$
{\p\vp\over\p x} = {\p (\vp_1,\ldots,\vp_n)\over \p(x_1,\ldots,x_n)} ,
$$
so that
$$
{\p\vp\over\p x} \cdot \left( {\p\vp\over\p x}\right)^\#
= \D\cdot I ,
$$
where $(\dscdot)^\#$ denotes the transposed matrix of
cofactors, $\D=\det(\p\vp/\p x)$ and $I$ is the identity
matrix.
Write $T_{(i)}^{(j)}$ for the $ji$ cofactor of $\p\vp/\p x$.
(Shrinking $V$ if necessary), we can choose constants
$c_1\ge 1$ and $c_2\ge 1$ such that, for all $\a\in\BN^n$ and
$a\in V$, $|D^\a \D (a)|\le \a! c_1 c_2^{|\a|}$ and each
$|D^\a T_{(i)}^{(j)} (a)| \le \a! c_1 c_2^{|\a|}$.

We can assume that $\vp(0)=0$.
Let $\mu(a)$ (respectively, $\nu(a)$) denote the order of
vanishing of $\D$ (respectively, of $(\p\vp/\p x)^\#$) at
a point $a\in V$.

Consider $f=g\circ\vp$, where $g$ is an analytic function in a
neighbourhood of $\vp(a)$, $a\in V$.
By the chain rule,
$$
\sum_{j=1}^n \left( {\p g\over\p y_j}\circ \vp\right)
{\p\vp_j\over \p x_i} = {\p f\over \p x_i} ,\qquad
i=1,\ldots,n .
$$
By Cramer's rule,
$$
\D \cdot \left( {\p\vp\over\p y_j}\circ \vp\right) =
\sum_{i=1}^n {\p f\over\p x_i} T_{(i)}^{(j)} ,\qquad
j=1,\ldots, n .
$$
 From the formulas obtained by repeated differentiation,
we get analytic functions $T_\a^\b$ on $V$, for all $\a,\b\in\BN^n$,
such that
$$
\D^{2|\b|-1}\cdot \big( (D^\b g)\circ\vp\big) =
\sum_{|\a|\le|\b|} T_\a^\b D^\a f ,
$$
where, at each $a\in V$, $T_\a^\b$ has order $\ge |\a|-\mu(a)+
|\b| \big( \mu(a)+\nu(a)-1\big)$ (cf. \cite{M, Lemme I.2.4}).
Moreover, there are constants $c_3\ge 1$, $c_4\ge 1$ depending
only on $c_1,c_2$ and $n$, such that, for all $\a,\b,\g\in \BN^n$
and all $a\in V$,
$$
\left| {D^\g T_\a^\b (a)\over \g!}\right| \le
c_3^{|\b|} c_4^{|\b|+|\g|-|\a|}
{(|\b|+|\g|-|\a|)!\over |\g|!}
$$
\cite{M, Lemme I.2.5}.

Let $a\in V$ and let $b=\vp(a)$.
If $F\in\hcO_a$ and $\a\in\BN^n$, let $D^\a F$ denote the
formal derivative of $F$ of order $\a$; thus, if $f\in\cO_a$ and
$\hf_a$ denotes the Taylor series of $f$ at $a$, then
$D^\a \hf_a=(D^a f)_a^{\widehat{\hphantom O}}$.
Consider $F=G\circ\hvp_a$, where $G$ is a formal power series
centred at $b$.
Write $G=\sum_{\b\in\BN^n} G_\b (y-b)^\b$, $F=\sum_{\a\in\BN^n}
F_\a (x-a)^\a$.
Assume that $F$ converges in a ball $|x-a|\le r$, where $r\le 1$;
thus, if $c=1/r$, then there is a constant $c'$ such that
$|F_\a|\le c' c^{|\a|}$, for all $\a\in\BN^n$.
As above, for all $\b\in\BN^n$, 
$$
\hDelta_a^{2|\b|-1} \cdot \big( (D^\b G)\circ \hvp_a\big) =
\sum_{|\a|\le |\b|} (T_\a^\b)_a^{\widehat{\hphantom O}}\cdot D^\a F .
\tag 3.1
$$
Let $H^\b$ denote the right-hand side of (3.1); write
$H^\b = \sum_{\g\in\BN^n} H_\g^\b (x-a)^\g$.
It is not difficult to estimate
$$
|H_\g^\b| \le 2^{|\b|+|\g|} (2n)^{|\b|} |\b|!
c_3^{|\b|} c_4^{|\b|+|\g|-1} c_5 c' c^{|\g|+\mu(a)-|\b|(\mu(a)+
\nu(a)-1)} ,
$$
where $c_5$ depends only on $c_3$, $c_4$ and $n$ \cite{M, I(2.6.9)}.

We now compare terms of degree $(2|\b|-1)\mu(a)$ in (3.1).
Order multiindices $\g\in\BN^n$ according to the lexicographic
order of $(|\g|,\g_1,\ldots,\g_n)$.
Let $\a$ denote the smallest $\g$ such that $D^\g \D(a)\ne 0$.
Then $(D^{(2|\b|-1)\a} \D^{2|\b|-1})(a) = \big( D^\a\D (a)\big)^
{2|\b|-1}$ and, from (3.1),
$$
\big( D^\a \D(a)\big)^{2|\b|-1} \b! G_\b = H_{(2|\b|-1)\a}^\b .
$$
We obtain
$$
|G_\b| \le c_6 c_7^{|\b|} c' c^{|\b|(\mu(a)-\nu(a)+1)} \cdot
{1 \over (D^\a \D(a))^{2|\b|-1}} ,
$$
where $c_6$, $c_7$ depend only on $n$, $\mu(a)$, $c_3$, $c_4$,
$c_5$.
In other words, $G$ converges in an open ball centred at
$b$, of radius at least
$$
{1\over c_7} |D^\a \D(a)|^2 r^{\mu(a)-\nu(a)+1} .
$$

Now, shrinking $V$ if necessary, we can assume that
$\mu(a)\le \mu(0)$ and $\nu(a)\le \nu(0)$, for all $a\in V$.
Therefore, $V$ has a finite partition into sets
$$
\Si_{\mu,\nu} := \big\{ a\in V:\ \mu(a)=\mu,\ \nu(a)=\nu
\big\} ;
$$
each $\Si_{\mu,\nu}$ is a difference of closed algebraic
(respectively, analytic) subsets of $V$.

For each $a\in V$, let $\a(a)$ denote the smallest multiindex
$\g$ such that $D^\g \D(a)\ne 0$ (in particular, $|\a(a)|=\mu(a)$).
Clearly, there is a finite filtration of $V$ by closed algebraic
(respectively, analytic) subsets $V=X_0\supset X_1\supset\cdots$
such that, for each $j=0,1,\ldots$, there exist $\mu$, $\nu$ and
$\a\in\BN^n$ such that $X_j\bs X_{j+1}\subset\Si_{\mu,\nu}$
and $\a(a)=\a$ for all $a\in X_j\bs X_{j+1}$.
Theorem 1.4 follows.

\enddocument